\documentclass[11pt]{amsart}
\usepackage{amsmath,amssymb,latexsym,dsfont}
\usepackage[left=2.6cm,right=2.6cm,top=2.7cm,bottom=2.7cm]{geometry}

\usepackage{amsmath,amsfonts,amssymb,epsfig, textcomp}
\usepackage{graphicx,tikz}
\usepackage{hyperref, cleveref} 
\usepackage{cases,listings}
\usepackage{color}
\usepackage{mathabx}

\newtheorem{theorem}{Theorem}[section]
\newtheorem{lemma}{Lemma}[section]

\newtheorem{proposition}{Proposition}[section]

\newtheorem{remark}{Remark}[section]

\numberwithin{equation}{section}

 \newcommand{\cB}{{\mathcal B}}

   \newcommand{\ty}{\widetilde y}

      \newcommand{\N}{\mathbb{N}}

      \newcommand{\eps}{\varepsilon}
      \newcommand{\mR}{\mathbb{R}}

      \newcommand{\mC}{\mathbb{C}}

      \newcommand{\supp}{\operatorname{supp}}

      \makeatletter
      \def\@setcopyright{}
      \def\serieslogo@{}
      \makeatother

\newcommand{\cG}{\mathcal G}

\newcommand{\tu}{\widetilde u}

\newcommand{\be}{\begin{equation}}
\newcommand{\ee}{\end{equation}}

\newcommand{\cF}{{\mathcal F}}

\newcommand{\cD}{{\mathcal D}}

\newcommand{\mH}{\mathbb{H}}

\title[Boundary control cost of one-dimensional fractional Schrödinger and heat equations]{Optimal bounds for the boundary control cost of one-dimensional fractional Schrödinger and heat equations}

\author[H.-M. Nguyen]{Hoai-Minh Nguyen}
\address[H.-M. Nguyen]{Sorbonne Universit\'e, Universit\'e Paris Cit\'e, CNRS, INRIA, \newline \indent 
Laboratoire Jacques-Louis Lions, LJLL, F-75005 Paris, France
}
\email{hoai-minh.nguyen@sorbonne-universite.fr}

\begin{document}

\maketitle 

\begin{center}
    \textit{Dedicated to Jean-Michel Coron on the occasion of his seventieth birthday, with admiration, gratitude, and friendship}
\end{center}

\begin{abstract} We derive sharp bounds for the boundary control cost of the one-dimensional fractional Schrödinger and heat equations. The analysis of the lower bound is based on the study of the control cost of a related singular boundary control problem in finite time, using tools from complex analysis. The analysis of the upper bound relies on the moment method, involving estimates of the Fourier transform of a class of compactly supported functions.
\end{abstract}

\bigskip 

\noindent {\bf Key words.} control cost, moment method, fractional heat equation, fractional Schr\"odinger equation. 

\noindent {\bf AMS subject classification.} 93B05, 93C20, 35B40.

\tableofcontents

\section{Introduction}

In this paper, we study the cost of the controls in small time of the fractional heat and fractional Schr\"odinger equation of order $1/ 2 < s < 1$ in one-dimensional space. More precisely, for $T> 0$, concerning the boundary controls for the Schr\"odinger equation, we deal with the system  
\be\label{sys-S-intro}
\left\{\begin{array}{cl}
i \partial_t y - (-\Delta)^s y = 0  \mbox{ in } (0, T) \times (0, 1), \\[6pt]
y(t, 0) = 0, \quad \cB y(t, 1) = u(t) \mbox{ in } (0, T), \\[6pt]
y(t=0, \cdot) = y_0 \mbox{ in } (0, 1),   
\end{array} \right.
\ee
and concerning the boundary control for the heat equation, we consider the control system 
\be \label{sys-H-intro}
\left\{\begin{array}{cl}
\partial_t y + (-\Delta)^s y = 0  \mbox{ in } (0, T) \times (0, 1), \\[6pt]
y(t, 0) = 0, \quad \cB y(t, 1) = u(t) \mbox{ in } (0, T), \\[6pt]
y(t=0, \cdot) = y_0 \mbox{ in } (0, 1).   
\end{array} \right.
\ee
Here, $y$ is the state, $y_0$ is the initial condition, and $u$ is the control. 

Some comments on these control systems are in order. 
Let $(\gamma_n)_{n \ge 1}$ be the eigenvalues of the $- \Delta$ operator with the zero Dirichlet boundary, and let $(\varphi_n)_{n \ge 1}$ be the corresponding eigenfunctions. We thus have, with $\Omega = (0, 1)$,  
$$
\left\{ \begin{array}{cl}
-\Delta \varphi_n = \gamma_n \varphi_n \mbox{ in } \Omega, \\[6pt]
\varphi_n = 0 \mbox{ on } \partial \Omega. 
\end{array} \right. 
$$ 
Assume that $(\varphi_n)_{n \ge 1}$ forms an orthogonal basis of $L^2(\Omega)$. Explicitly, one has
$$
\gamma_n = \pi^2 n^2  \mbox{ for } n \ge 1, 
$$
and one can choose, for $n \ge 1$,  
$$
\varphi_n = \sqrt{2} \sin (\pi n x) \mbox{ for } x \in (0, 1). 
$$

Denote $\mH = L^2(\Omega)$ and $\langle \cdot, \cdot \rangle$ the standard scalar product in $\mH$. Define
$$
\cD( (-\Delta)^s) = \left\{\varphi \in \mH;  \sum_{n \ge 1} \gamma_n^{2s} |\langle \varphi, \varphi_n \rangle |^2 < + \infty \right\}
$$
and
$$
(-\Delta)^s \varphi =   \sum_{n \ge 1} \gamma_n^{s} \langle \varphi, \varphi_n \rangle \varphi_n \mbox{ for } \varphi \in \cD( (-\Delta)^s). 
$$
The control system \eqref{sys-S-intro} is understood under the sense that $y \in C([0, T]; L^2(0, 1))$, and $y$ satisfies 
\be \label{thm-S-meaning}
i \frac{d}{dt} \langle y, \varphi \rangle - \langle y, (-\Delta)^s \varphi \rangle = u(t) \varphi_x(1) \mbox{ in } (0, T) \mbox{ for $\varphi = \varphi_n$ with $n \ge 1$}, 
\ee
in the distributional sense, and the 
control system \eqref{sys-H-intro} is understood under the sense that $y \in C([0, T]; L^2(0, 1))$, and $y$ satisfies  
\be \label{thm-H-meaning}
\frac{d}{dt} \langle y, \varphi \rangle + \langle y, (-\Delta)^s \varphi \rangle =-  u(t)  \varphi_x(1) \mbox{ in } (0, T) \mbox{ for $\varphi = \varphi_n$ with $n \ge 1$}, 
\ee
in the distributional sense. When $s=1$, we rediscover the boundary control of Schr\"odinger equation using the Dirichlet control on the right.

\medskip

Concerning the fractional Schr\"odinger system \eqref{sys-S-intro}, one has the following result.  

\begin{theorem} \label{thm-S} Let $1/2 < s  < 1$ and $0< T < 1/2$, and let  $y_0 \in L^2(0, 1)$. Set
\be
\tau = \frac{1}{2 s - 1}. 
\ee
We have
\begin{itemize}
\item[a)] Assume that $y_0 = \varphi_1$ and  $u \in L^2(0, T)$ is such that $y(T, \cdot) = 0$ where 
$y \in C([0, T]; L^2(0, 1))$ is the unique solution of the system \eqref{sys-H-intro}. Then 
\be
\| u\|_{L^2(0, T)} \ge T^{c_1} e^{\frac{\nu_s}{T^{\tau}}} \| y_0\|_{L^2(0, T)}. 
\ee
\item[b)] There exists $u \in L^2(0, T)$ such that $y(T, \cdot) = 0$ where 
$y \in C([0, T]; L^2(0, 1))$ is the unique solution of the system \eqref{sys-S-intro}. Moreover, one can choose $u$ such that 
\be
\| u\|_{L^2(0, T)} \le c_2 e^{\frac{c_3}{T^{\tau}}} \| y_0\|_{L^2(0, T)}. 
\ee
\end{itemize}
Here $c_1, c_2, c_3$ are positive constants depending only on $s$, and the constant $\nu_s$ is defined by 
\be
\nu_s = \frac{1}{2} (2s-1)\left( \frac{1}{2 s\sin (\pi / (4s) )} \right)^{\frac{2s}{2s-1}}.
\ee
\end{theorem}

One can check that 
\be
\nu_s =  \left( \frac{\theta_{2s}}{\pi} \right)^{1 + 1/ \beta} \left(\frac{2}{\beta + 1}\right)^{1/\beta}  \frac{\beta}{1 + \beta} \mbox{ with } \beta = 2 s - 1, 
\ee
where, for $\alpha > 1$,  
\be \label{def-cs}
\theta_{\alpha} = - \int_0^\infty \ln \left( \frac{x^{\alpha}}{\sqrt{x^{2 \alpha} + 1}} \right) \, dx = \frac{1}{2} \int_0^\infty \ln (1 + x^{-2\alpha}) \, dx = \frac{\pi}{2 \sin (\pi/ (2\alpha) )} . 
\ee

Concerning the fractional heat system \eqref{sys-H-intro}, one has the following result.  

\begin{theorem} \label{thm-H} Let $1/2 < s < 1$ and $0< T < 1/2$, and let  $y_0 \in L^2(0, 1)$. Set
$$
\tau = \frac{1}{2 s - 1}. 
$$
We have
\begin{itemize}
\item[a)] Assume that $y_0 = \varphi_1$ and  $u \in L^2(0, T)$ is such that $y(T, \cdot) = 0$ where 
$y \in C([0, T]; L^2(0, 1))$ is the unique solution of the system \eqref{sys-H-intro}. Then 
\be
\| u\|_{L^2(0, T)} \ge T^{c_1} e^{\frac{\mu_s}{T^{\tau}}} \| y_0\|_{L^2(0, T)}. 
\ee
\item[b)] There exists $u \in L^2(0, T)$ such that $y(T, \cdot) = 0$ where 
$y \in C([0, T]; L^2(0, 1))$ is the unique solution of the system \eqref{sys-H-intro}. Moreover, one can choose $u$ such that 
\be
\| u\|_{L^2(0, T)} \le c_2 e^{\frac{c_3}{T^{\tau}}} \| y_0\|_{L^2(0, T)}. 
\ee
\end{itemize}
Here $c_1, c_2, c_3$ are positive constants depending only on $s$, and the constant $\mu_s$ is defined by 
\be
\mu_s =  \frac{1}{2} (2s-1)\left( \frac{1}{s\sin (\pi / (2s) )} \right)^{\frac{2s}{2s-1}}. 
\ee
\end{theorem}

One can check that 
\be
\mu_s =  \left( \frac{\kappa_{2s}}{\pi} \right)^{1 + 1/ \beta} \left(\frac{2}{\beta + 1}\right)^{1/\beta}  \frac{\beta}{1 + \beta} \mbox{ with } \beta = 2s -1,  
\ee
where, for $\alpha > 1$,  
\be \label{def-ds}
\kappa_\alpha = - \int_0^\infty \ln \left( \frac{x^{\alpha}}{x^{\alpha} + 1} \right) \, dx = \frac{\pi}{ \sin (\pi / \alpha)}. 
\ee

These results can be proved under more general assumptions on  $(\lambda_n)_{n \ge 1}$ where $\lambda_n = \gamma_n^{\alpha/2}$ with $\alpha = 2s$. Indeed, these hold for the following assumption on $(\lambda_n)$: 
\be \label{assumption-lambda_n}
\gamma = \inf_{k \neq n} |\lambda_k - \lambda_n|, \quad  \Gamma_1 = \sup_{k} \frac{|\lambda_k - a k^\alpha|}{k^{\alpha - 1}}, \quad \mbox{ and } \quad \Gamma_2 = \sup \frac{k^\alpha}{\lambda_k}.   
\ee
Moreover, one can show that the constants appearing here can be chosen in such a way that they depend only on $\alpha$, $a$, $\gamma$, $\Gamma_1$, and $\Gamma_2$ (see \Cref{lem-S,lem-H,lem-Phi,lem-HH}). 

The null-controllability of the Schr\"odinger equation ($s = 1$) was established by  Lebeau \cite{Lebeau92} under the geometric control condition. This geometric control condition involving geometric optics was first considered for the wave equation of Rauch and Taylor \cite{RT74} and Bardos, Lebeau, and Rauch \cite{BLR92}. Variants of \Cref{thm-S} are known in the case $s = 1$, see, e.g., \cite{Miller04,Miller04-JDE,TT07} and the references therein. It was shown by Biccari \cite{Biccari22} using Ingham’s inequalities that the fractional Schr\"odinger equation in one dimension is controllable for $1/2 < s < 1$ in small time using internal controls.

The null-controllability of the heat equation ($s = 1$) is well-known. In one dimensional space, it was obtained  by Fattorini and Russell \cite{FR71} where the moment method was introduced. In higher dimensions, the null-controllability was established by Lebeau and Robiano \cite{LR95} via spectral inequalities and Fursikov and Imanuvilov \cite{FI96} via Carleman estimates. Variants of \Cref{thm-H} are known in the case $s = 1$, see, e.g., \cite{Guichal85,SAI00,Miller04,Miller04-JDE,TT07}.  In the case $1/2< s < 1$, Micu and Zuazua \cite{MZ06} using the moment method proved that the system \eqref{sys-H-intro} is null-controllable for any time $T>0$ \footnote{In fact, they considered lump controls (bilinear controls), but their method and results also hold in this case.}. They also proved that this system is not null-controllable for $s < 1/2$. Concerning the internal controls with $1/2 < s < 1$, Miller \cite{Miller06} established that the system is null-controllable and the control cost is bounded  $c_1 e^{c_2/T^{\tau'}}$ for $\tau' > \tau$ in arbitrary dimensions using Lebeau \& Robiano's approach (see also \cite{MZ06}). Koenig \cite{Koenig20} showed that the system is not null-controllable for $s < 1/2$ for a torus in arbitrary dimensions (see also \cite{Koenig17} for the case $s=1/2$ in one dimensional case). 

Part $a)$ of \Cref{thm-S,thm-H} was obtained by Lissy \cite{Lissy15}. Part $b)$ of \Cref{thm-S,thm-H} was obtained by Lissy \cite{Lissy17} under assumption \eqref{assumption-lambda_n}.   In the case $s > 1$, similar results  were established by Lissy \cite{Lissy14} (see also \cite{Lissy15} for the lower bound and \cite{Lissy17} for the upper case).

The proof of the lower bound is inspired by the work of Coron and Guerrero \cite{CG05} (see also \cite{Coron07}). The idea is to link the cost for fast control of the control problems considered with singular ones for which the tools from complex analysis, see, e.g., \cite{Koosis09}, 
can be implemented as previously suggested by Coron and Guerrero \cite{CG05}.  Indeed, we first study the lower bounds of the control cost of the following two systems, with $\beta = 2 s -1$, 
\be \label{sys-S-singular}
i \partial_t y - \eps^\beta (-\Delta)^s y  + \frac{i }{2 \eps } y = 0  \mbox{ in } (0, T) \times (0, 1), 
\ee
(see \Cref{pro-S-M}) and 
\be \label{sys-H-singular}
y_t +  \eps^\beta (-\Delta)^s y + \frac{1}{2 \eps } y = 0   \mbox{ in } (0, T) \times (0, 1), 
\ee
(see \Cref{pro-H-M}) instead of 
$$
i \partial_t y - (-\Delta)^s y = 0  \mbox{ in } (0, T) \times (0, 1), 
$$
and 
$$
\partial_t y + (-\Delta)^s y = 0  \mbox{ in } (0, T) \times (0, 1), 
$$
respectively. We then derive the lower bound for the cost of controls of the original systems. The method given here also works for the case $s \ge 1$ (see \Cref{lem-S,lem-HH}). The analysis given here share many common points with one of Lissy \cite{Lissy15} where he considered the case $\lambda_n = a n^\alpha$. Thus the lower bound results under the assumption \eqref{assumption-lambda_n} is new to our knowledge.  An effort given  to derive the same estimate under the assumption \eqref{assumption-lambda_n} where the dependence is on $\alpha$, $a$, $\gamma$, $\Gamma_1$, and $\Gamma_2$. 

The proof of the upper bound is based on the moment method. This method has been used to study the controllability and the cost of the controls in various control problems in one dimension, see, e.g.,  \cite{CG05,MZ06,TT07,BL10,G10,Lissy14,ABGT14} and the references therein. Our analysis has its roots in the work of Tenenbaum and Tucsnak \cite{TT07}. To handle the fractional case $1/2 < s < 1$, two additional tasks in comparison with \cite{TT07} are required.  
First, we need to construct explicitly a function $H$ with support in $[-T/2, T/2]$ whose Fourier transforms decays at infinity  as $e^{-c |x|^{\gamma}}$ for $1/2 < \gamma < 1$ (see \Cref{lem-H}). The case $\gamma = 1/2$ was previously considered in \cite{TT07}. Second, we need to estimate the behavior of holomorphic functions $\Phi$, whose zeros are given by $\lambda_k: = \gamma_k^s$ or $i \lambda_k: = i \gamma_k^s$ for $k \ge 1$. These are based on \Cref{lem-Phi} which requires an estimate on the infinite product given in
\Cref{lem-product} with the roof from the work of Fattorini and Russell \cite{FR71}. 
Our arguments are simpler and different from \cite{TT07}. It is worth mentioning that the direct uses of the results and arguments in \cite{TT07} cannot handle the case $1/2< s < 1$ as previously discussed in \cite{Lissy14} (see \cite[Remarks 3 and 4]{Lissy14} and see also \cite[Lemmas 2.1 and 2.3]{Lissy14}).  Our proof is thus different with that of Lissy given in \cite{Lissy17}, even though both are based on the moment method. 
The construction of $H$ in \cite{Lissy17} is based on the Bray-Mandelbrojt construction of Gevrey functions and is somehow more complicated than ours. Our estimate of $\Phi$ is more elementary than that of \cite{Lissy17} and the details on the  information of the infinite product are now given (see \Cref{lem-product}). We do not concentrate on optimizing parameters to obtain the best possible constants for the upper bound since the method used cannot give a better constant than some other methods, see the recent works of Dard\'e and Ervedoza \cite{DE19} on the heat equation, and the work of Lissy \cite{Lissy25} on the Schr\"odinger equation.

Good understanding the cost of fast controls for control systems generated by strongly continuous {\it group} can lead to a process to stabilize the system in finite time using Gramian operators. This way has been introduced and implemented for the KdV equations using boundary controls in \cite{Ng-S-KdV} and the Schr\"odinger equation using bilinear controls in \cite{Ng-S-Schrodinger}, with the roots from \cite{Ng-Riccati}. It is worth noting that, the fractional heat equation is not null-controllable in the case $0 < s < 1/2$. Nevertheless, one can still achieve the rapid stabilization in this case using backstepping technique, see \cite{GKN25}.

The paper is organized as follows. The lower bounds (part $a)$) given in \Cref{thm-S} and \Cref{thm-H} are proved in \Cref{sec-LB-S} and \Cref{sec-LB-H}, respectively.  In \Cref{sec-lemmas}, we present some technical lemmas on the functions $\Phi$ and $H$ used in the proof of the upper bound. The upper bounds (part $b)$) given in \Cref{thm-S} and \Cref{thm-H} are proved in \Cref{sec-upperbound}.

\section{Lower bound on the cost of fast controls for the Schrödinger equation} \label{sec-LB-S}

\subsection{Preliminaries}

In this section, we established two lemmas. Both of them can be used to derive part $a)$ of \Cref{thm-S}. The first lemma, can be also used to deal with \eqref{sys-S-singular}, is the following. 

\begin{lemma} \label{lem-S-M} Let $\alpha > 1$, $a > 0$,  and let $(\lambda_k)_{k \ge 1}$ be a positive sequence. 
Assume that $\gamma > 0$ and $\Gamma_1 < + \infty$, where 
\be \label{lem-S-M-assumption}
\gamma = \mathop{\inf_{k \neq n}}_{k, n \ge 1} |\lambda_k - \lambda_n| \quad \mbox{ and } \quad \Gamma_1 = \sup_{k \ge 1} \frac{|\lambda_k - a k^\alpha|}{k^{\alpha - 1}}.  
\ee
Let $T>0$ and  $u \in L^1(0, T)$. Define 
\be
\cG (z) = \int_0^T u(t) e^{- i z t} \, dt. 
\ee
Denote 
$$
\beta = \alpha - 1, 
$$
and set, for $0< \eps < m$,  
\begin{equation}
\Lambda_k :=  - \eps^\beta \lambda_k +  \frac{i}{2 \eps}.
\end{equation}
Assume that 
\be
\cG(\Lambda_1) = 1 \quad \mbox{ and } \quad \cG (\Lambda_k) = 0 \mbox{ for } k \ge 2.
\ee 
Then 
\begin{equation}
\| u \|_{L^1(0, T)}   \ge  c \eps^C e^{\frac{1}{\eps} \left(\frac{\theta_\alpha}{a^{1/\alpha}} - \frac{T}{2} \right)} 
\end{equation}
for some positive constants $c, C$ depending only $\alpha$, $a$, $\gamma$, $\Gamma_1$, and $m$. Here $\theta_\alpha$ is defined by \eqref{def-cs}.  

\end{lemma}

\begin{proof} We follow the analysis suggested by Coron and Guerrero \cite{CG05}. 
Define, for $z \in \mC$,
\begin{equation}\label{pro-S-M-F}
\cF(z) :=   \int_{-T/2}^{T/2} v(t) e^{- i z t} \, dt \mbox{ where } v(t) = u(t + T/2) \mbox{ for } t \in (-T/2, T/2). 
\end{equation}
It follows that 
\begin{equation}\label{pro-S-M-zeros}
\cF (\Lambda_{1}) = e^{ i \Lambda_1 T/2}
\quad \mbox{ and } \quad \cF (\Lambda_k ) =0  \mbox{ for } k \neq 1.
\end{equation}
Applying  the representation of entire functions of exponential type for $\cF$, see e.g.,  \cite[page 56]{Koosis09}, we derive from \eqref{pro-S-M-zeros} that, for $z  \in \mC$ with $\Im z > 0$,
\begin{equation}\label{pro-S-M-complex-analysis} 
\ln |\cF (z)| \le I_0(z)  + I_1(z) + \sigma \Im( z),
\end{equation}
where
\begin{equation} \label{pro-S-M-I01}
I_0 (z)= \sum_{k \neq 1} \ln \frac{\big|\Lambda_k - z \big|}{\big|
 \overline{\Lambda_k} - z \big|}, \quad I_1(z) = \frac{\Im(z)}{\pi} \int_{-\infty}^\infty \frac{\ln |\cF(\tau)|}{|\tau - z|^2} \, d\tau,
\end{equation}
and
\begin{equation}
\sigma = \limsup_{y \to + \infty}  \frac{\ln |\cF( i y)|}{y}.
\end{equation}

From the definition of $\cF$ in \eqref{pro-S-M-F}, we have, for $y \in \mR_+$,  
\begin{equation*}
|\cF(i y)| \le \| u\|_{L^1(0, T)} e^{T y/2}. 
\end{equation*}
This implies 
\begin{equation}\label{pro-S-M-est1}
\sigma = \limsup_{y \to + \infty}  \frac{\ln |\cF(i y)|}{y} \le T/2. 
\end{equation}

From  \eqref{pro-S-M-zeros},  we derive that 
\begin{equation}\label{pro-S-M-lnF}
\ln |F(\Lambda_{1})| 
\ge  - \frac{T}{4 \eps}.
\end{equation}

By \eqref{pro-S-M-I01}, we have 
\be\label{pro-S-M-est-I0-1}
I_0 (\Lambda_{1} ) 
\le \sum_{k \ge 2} \ln \left( \frac{\eps^\beta |\lambda_k -\lambda_1| }{ \sqrt{(\eps^\beta (\lambda_k - \lambda_1))^2 + (1/\eps)^2}} \right) = \sum_{k \ge 2} \ln \left( \frac{ \eps^{\beta+1} |\lambda_k -\lambda_1|  }{ \sqrt{(\eps^{\beta+1} (\lambda_k - \lambda_1))^2  + 1}} \right).  
\ee
In what follows,  $C$ denotes a positive constant depending only on $\alpha$, $a$, $\gamma$, $\Gamma_1$, and $m$, and can change from one place to another. 
It follows from \eqref{lem-S-M-assumption} that
\begin{multline}\label{pro-S-M-est-I0-2}
I_0 (\Lambda_{1} ) \le  \sum_{k \ge C} \ln \left( \frac{ a \eps^{\beta+1} k^\alpha }{ \sqrt{(a \eps^{\beta+1} k^\alpha)^2  + 1}} \right)
\\[6pt] 
\le   \int_{C}^\infty \ln \left(\frac{a \eps^{\beta +1}  x^\alpha}{\sqrt{a^2 \eps^{2(\beta + 1)}  x^{2\alpha} + 1 }} \right) \, dx 
\le \frac{1}{a^{\frac{1}{\alpha}} \eps^{\frac{\beta + 1}{\alpha}}} \int_{C a^{1/\alpha} \eps^{\frac{\beta + 1}{\alpha}}}^\infty \ln \left(\frac{ x^\alpha}{\sqrt{x^{2\alpha} + 1}} \right) \, dx. 
\end{multline}
 Since 
$$
\int_{0}^\infty \ln \Big(\frac{ x^\alpha}{\sqrt{x^{2\alpha} + 1}} \Big)  \, dx = -  \theta_\alpha, 
$$
it follows from \eqref{pro-S-M-est-I0-2} that
\begin{equation}\label{pro-S-M-I0}
I_0 (\Lambda_{1} )  \le - \frac{\theta_\alpha}{a^{\frac{1}{\alpha}} \eps^{\frac{\beta + 1}{\alpha}}} - C \ln \eps + C = - \frac{\theta_\alpha}{a^{\frac{1}{\alpha}} \eps} - C \ln \eps + C. 
\end{equation}

From \eqref{pro-S-M-F},  we have, for $s \in \mR$, 
$$
\ln |\cF(s)|  \le \| u\|_{L^1(0, T)},   
$$
which yields 
\begin{equation}\label{pro-S-M-I1-1}
I_1(\Lambda_{1}) \le \frac{1}{2 \eps \pi} \int_{-\infty}^\infty \frac{\ln \|u \|_{L^1(0, T)} }{s^2 + 1/(2\eps)^2} \, ds. 
\end{equation}
Since, for $\xi>0$,
\begin{equation*}
\frac{\xi}{\pi} \int_{-\infty}^\infty \frac{1}{s^2 + \xi^2} \, ds = 1,
\end{equation*}
it follows from \eqref{pro-S-M-I1-1} that
\begin{equation}\label{pro-S-M-I1}
I_1 \big(\Lambda_{1} \big) \le \ln  \|u \|_{L^1(0, T)}.
\end{equation}

Combining \eqref{pro-S-M-complex-analysis}, \eqref{pro-S-M-lnF}, \eqref{pro-S-M-I0}, and \eqref{pro-S-M-I1}  yields
\begin{equation*}
  \frac{T}{4\eps} \le - \frac{\theta_\alpha}{a^{\frac{1}{\alpha}} \eps}  + \ln \| u \|_{L^1(0, T)}  + \frac{T}{4 \eps} - C \ln \eps + C.  
\end{equation*}
This implies
\begin{equation*}
\ln  \| u \|_{L^1(0, T)}   \ge \frac{\theta_\alpha}{a^{\frac{1}{\alpha}} \eps}  - \frac{T}{2 \eps}  + C \ln \eps - C. 
\end{equation*}
The conclusion follows.  
\end{proof}

Here is the second lemma, which is a consequence of  \Cref{lem-S-M}. 

\begin{lemma} \label{lem-S} Let $\alpha > 1$, $a > 0$,  and let $(\lambda_k)_{k \ge 1}$ be a positive sequence. 
Assume that $\gamma > 0$ and $\Gamma_1 < + \infty$, where 
\be
\gamma = \inf_{k \neq n} |\lambda_k - \lambda_n| \quad \mbox{ and } \quad \Gamma_1 = \sup_{k} \frac{|\lambda_k - a k^\alpha|}{k^{\alpha - 1}}. 
\ee
Let $0< T < m$ and  $u \in L^1(0, T)$. Define 
\be
\cG (z) = \int_0^T u(t) e^{- i z t} \, dt. 
\ee
Assume that $\cG(- \lambda_1) = 1$ and $\cG (-\lambda_k) = 0$ for $k \ge 2$. Then 
\begin{equation}
\| u \|_{L^1(0, T)}   \ge  c T^{C} e^{\frac{\rho}{T^{\tau}}} \mbox{ with } \tau = \frac{1}{\alpha - 1}, 
\end{equation}
for some positive constants $c$, $C$ depending only on $\alpha, \, a,  \,  \gamma, \,  \Gamma_1$, and $m$. Here 
$$
\rho = 2^{\frac{1}{\beta}} \beta \left(  \frac{\theta_\alpha}{a^{1/\alpha} (1 + \beta)}\right)^{1 + \beta}  \mbox{ with } \beta = \alpha - 1,  
$$
with $\theta_\alpha$ being defined by \eqref{def-cs}.  

\end{lemma}

\begin{proof} Denote 
$$
T_0 = \frac{2 \theta_\alpha}{a^{\frac{1}{\alpha}} (\beta + 1)}, \quad \eps^\beta T_0 = T, 
$$
and set 
\be
\tu(t) = u(\eps^\beta t) e^{-\frac{t}{2 \eps}} \mbox{ for } t \in (0, T_0). 
\ee
Define, for $z \in \mC$,  
\be
\widetilde G(z) = \int_0^{T_0} \tu (t) e^{- i z t} \, dt. 
\ee
We have 
\begin{multline}
G(z) = \int_0^T u(t) e^{- i z t} \, dt =  \int_0^{T} \tu (\eps^{-\beta} t) e^{\eps^{- \beta} t/ (2 \eps)} e^{- i z t} \, dt \\[6pt]
 = \int_0^{T} \tu (\eps^{-\beta} t) e^{\eps^{- \beta} t/ (2 \eps)} e^{- i \eps^\beta z (\eps^{-\beta } t) } \, dt 
= \eps^{\beta} \int_0^{\eps^{-\beta} T} \tu (t)   e^{  t/ (2 \eps)} e^{- i \eps^\beta z  t } \, dt \\[6pt]
= \eps^{\beta}   \int_0^{T_0} \tu (t)  e^{- i \Big(\eps^\beta z + i / (2 \eps) \Big)  t } \, dt = \eps^\beta \widetilde G (\eps^\beta z + i / (2 \eps)). 
\end{multline}
Applying \Cref{lem-S-M}, we have 
\be
\| \tu\|_{L^1(0, T_0)} \ge c \eps^{C} e^{\frac{1}{\eps} \left( \frac{\theta_\alpha}{a^{1/ \alpha}} - \frac{T_0}{2}\right)}. 
\ee 
Since 
$$
\frac{1}{\eps} \left( \frac{\theta_\alpha}{a^{1/ \alpha}} - \frac{T_0}{2}\right) = \frac{\rho}{T^\tau},  
$$
the conclusion follows. 
\end{proof}

\subsection{Lower bound on the cost of fast controls for a singular Schr\"odinger control system}

In this section, we deal with \eqref{sys-S-singular}.

\begin{proof} Define, for $k \ge 1$, 
\begin{equation}
\Lambda_k :=  - \eps^\beta \lambda_k +  \frac{i}{2 \eps}, 
\end{equation}
where
\be
\lambda_k = \gamma_k^{\alpha/2} \mbox{ with } \alpha = 2 s.  
\ee

Consider \eqref{pro-S-M-meaning} with  $\varphi = \varphi_k$. Multiplying this system with $e^{- i \Lambda_k t}$ and integrating over $(0, T)$, we obtain
\be
- i  \langle y_0, \varphi_k \rangle = (-1)^k k \pi \int_0^T u(t) e^{- i \Lambda_k t} \, dt.
\ee

Define, for $z \in \mC$,
\begin{equation}
\cG(z) =   \int_{0}^{T} u(t) e^{- i z t} \, dt. 
\end{equation}
It follows that 
\begin{equation}
\cG (\Lambda_{1}) = \frac{1}{i \pi}
\quad \mbox{ and } \quad \cG (\Lambda_k ) =0  \mbox{ for } k \ge 2.
\end{equation}
The conclusion now follows from \Cref{lem-S-M}.  
\end{proof}

\subsection{Proof of part $a)$ of \Cref{thm-S}} Denote 
$$
\beta = 2s  - 1, \quad T_0 = \frac{2 \theta_{2s}}{\pi (\beta + 1)}, \quad \eps^\beta T_0 = T, 
$$
and set 
\be
\ty(t, x) = y(\eps^\beta t, x) e^{-\frac{t}{2 \eps}},  \mbox{ and } \quad \tu(t) = u(\eps^\beta t) e^{-\frac{t}{2 \eps}}. 
\ee
We then have 
\begin{equation} \left\{
\begin{array}{c}
i \ty_t -  \eps^\beta (-\Delta)^s \ty + \frac{i}{2 \eps } \ty = 0   \mbox{ in } (0, T_0) \times (0, L), \\[6pt]
\ty(t, 0) = 0 \mbox{ in } (0, T_0), \\[6pt]
\ty (t, L) = \tu(t) \mbox{ in } (0, T_0), \\[6pt]
\ty_0 = \varphi_{1} \mbox{ in } (0, 1). 
\end{array} \right. 
\end{equation} 
Moreover, 
$$
\ty(T_0, \cdot) = y (T, \cdot) e^{-\frac{T_0}{2 \eps}} = 0. 
$$
By \Cref{pro-S-M}, we have 
$$
\|\tu \|_{L^1(0, T_0)} \ge \eps^C e^{\frac{1}{\eps} \left(\frac{\theta_{2s}}{\pi} - \frac{T_0}{2} \right)}. 
$$
Since 
\be
\frac{1}{\eps} \left(\frac{\theta_{2s}}{\pi} - \frac{T_0}{2} \right)= \frac{\nu_s}{T^\tau}, 
\ee
the conclusion follows. \qed

\begin{remark} \rm
One can also obtain a proof of part $a)$ of \Cref{thm-S} using \Cref{lem-S}. 
\end{remark}

\section{Lower bound on the cost of fast controls for the heat equation}\label{sec-LB-H}

\subsection{Preliminaries}

In this section, we established two lemmas. Both of them can be used to derive part $a)$ of \Cref{thm-H}. The first lemma, can be used to deal with \eqref{sys-H-singular} is the following.

\begin{lemma} \label{lem-H-M} Let $\alpha > 1$, $a > 0$,  and let $(\lambda_k)_{k \ge 1}$ be a positive sequence. 
Assume that $\gamma > 0$ and $\Gamma_1 < + \infty$ where 
\be \label{lem-H-M-assumption}
\gamma = \inf_{k \neq n} |\lambda_k - \lambda_n| \quad \mbox{ and } \quad \Gamma_1 = \sup_{k} \frac{|\lambda_k - a k^\alpha|}{k^{\alpha - 1}}.  
\ee
Let $T>0$ and  $u \in L^1(0, T)$. Define, for $z \in \mC$,  
\be
\cG (z) = \int_0^T u(t) e^{- i z t} \, dt. 
\ee
Denote 
$$
\beta = \alpha - 1
$$
and set, for $0< \eps < m$,  
\begin{equation*}
\Gamma_k :=   \eps^\beta \lambda_k +  \frac{1}{2 \eps}, 
\end{equation*}
Assume that $\cG(i \Gamma_1) = 1$ and $\cG (i \Gamma_k) = 0$ for $k \ge 2$. Then 
\begin{equation}
\| u \|_{L^1(0, T)}   \ge  c \eps^C e^{\frac{1}{\eps} \left(\frac{\kappa_\alpha}{a^{1/\alpha}} - \frac{T}{2} \right)} 
\end{equation}
for some positive constants $c$,  $C$ depending only on $\alpha, \, a, \,  \gamma, \, \Gamma_1$ and $m$. Here $\kappa_\alpha$ is defined by \eqref{def-ds}.  
\end{lemma}

\begin{proof} The proof is in the same spirit of the one of \Cref{lem-S-M}. Define, for $z \in \mC$,
\begin{equation}\label{pro-H-M-F}
\cF(z) :=   \int_{-T/2}^{T/2} v(t) e^{- i z t} \, dt \mbox{ with } v(t) = u(t + T/2) \mbox{ for } t \in (-T/2, T/2). 
\end{equation}
It follows that 
\begin{equation}\label{pro-H-M-zeros}
\cF (i \Gamma_{1}) = e^{\Gamma_1 T/2}
\quad \mbox{ and } \quad \cF (i \Gamma_k ) =0  \mbox{ for } k \neq 2.
\end{equation}

Applying  the representation of entire functions of exponential type for $\cF$, see e.g.,  \cite[page 56]{Koosis09}, we derive from \eqref{pro-H-M-zeros} that, for $z  \in \mC$ with $\Im z > 0$,
\begin{equation}\label{pro-H-M-complex-analysis} 
\ln |\cF (z)| \le I_0(z)  + I_1(z) + \sigma \Im( z),
\end{equation}
where
\begin{equation} \label{pro-H-M-I01}
I_0 (z)= \sum_{k \neq 1} \ln \frac{\big|i \Gamma_k - z \big|}{\big|
- i  \Gamma_k - z \big|}, \quad I_1(z) = \frac{\Im(z)}{\pi} \int_{-\infty}^\infty \frac{\ln |\cF(t)|}{|t - z|^2} \, dt,
\end{equation}
and
\begin{equation}
\sigma = \limsup_{y \to + \infty}  \frac{\ln |\cF( i y)|}{y}.
\end{equation}

From the definition of $\cF$ in \eqref{pro-H-M-F}, we have, for $y \in \mR_+$,  
\begin{equation*}
|\cF(i y)| \le \| u\|_{L^1(0, T)} e^{T y/2}. 
\end{equation*}
This implies 
\begin{equation}\label{pro-H-M-est1}
\sigma = \limsup_{y \to + \infty}  \frac{\ln |\cF(i y)|}{y} \le T/2. 
\end{equation}

From  \eqref{pro-H-M-zeros},  we derive that 
\begin{equation}\label{pro-H-M-lnF}
\ln |F(i \Gamma_{1})| 
\ge  - \frac{T}{4 \eps} - C.
\end{equation}

By \eqref{pro-H-M-I01}, we have 
\be
I_0 (i \Gamma_{1} ) 
\le  \sum_{k \ge 2} \ln \frac{\eps^\beta |\lambda_k - \lambda_1| }{ \eps^\beta |\lambda_k + \lambda_1| + 1/\eps} = \sum_{k \ge 2} \ln \frac{ \eps^{\beta+1} |\lambda_k - \lambda_1| }{ \eps^{\beta+1} |\lambda_k + \lambda_1|  + 1}. 
\ee
In what follows,  $C$ denotes a positive constant depending only on $\alpha$, $a$, $\gamma$, $\Gamma_1$, and $A$, and can change from one place to another. 
It follows from \eqref{lem-H-M-assumption} that
\begin{multline}\label{pro-H-M-est-I0-1}
I_0 (i \Gamma_{1} )  \le  \sum_{k \ge C} \ln \frac{\eps^\beta a k^\alpha }{ \eps^\beta a k^\alpha+ 1/\eps}
\le   \int_{C}^\infty \ln \left(\frac{ \eps^{\beta +1} a x^\alpha}{ \eps^{\beta + 1} a x^{\alpha} + 1 } \right) \, dx 
\le \frac{1}{a^{\frac{1}{\alpha}} \eps^{\frac{\beta + 1}{\alpha}}} \int_{Ca^{1/\alpha} \eps^{\frac{\beta + 1}{\alpha}}}^\infty \ln \Big(\frac{ x^\alpha}{x^\alpha + 1} \Big) \, dx. 
\end{multline}
Since 
$$
\int_{0}^\infty \ln \Big(\frac{x^\alpha}{x^\alpha + 1} \Big) \, dx = -  \kappa_\alpha, 
$$
it follows from \eqref{pro-H-M-est-I0-1} that
\begin{equation}\label{pro-H-M-I0}
I_0 (i \Gamma_{1} )  \le - \frac{\kappa_\alpha}{a^{\frac{1}{\alpha}} \eps^{\frac{\beta + 1}{\alpha}}} - C \ln \eps  + C = - \frac{\kappa_\alpha}{a^{\frac{1}{\alpha}} \eps} - C \ln \eps + C. 
\end{equation}

From \eqref{pro-H-M-F},  we have, for $t \in \mR$, 
$$
\ln |\cF(t)|  \le \| u\|_{L^1(0, T)},   
$$
which yields 
\begin{equation}\label{pro-H-M-I1-1}
I_1(i\Gamma_{1}) \le \frac{1}{2 \eps \pi} \int_{-\infty}^\infty \frac{\ln \|u \|_{L^1(0, T)} }{t^2 + 1/(2\eps)^2} \, dst. 
\end{equation}
Since, for $a>0$,
\begin{equation*}
\frac{a}{\pi} \int_{-\infty}^\infty \frac{1}{t^2 +a^2} \, dt = 1,
\end{equation*}
it follows from \eqref{pro-H-M-I1-1} that
\begin{equation}\label{pro-H-M-I1}
I_1 \big(i \Gamma_{1} \big) \le \ln  \|u \|_{L^1(0, T)}.
\end{equation}

Combining \eqref{pro-H-M-complex-analysis}, \eqref{pro-H-M-lnF}, \eqref{pro-H-M-I0}, and \eqref{pro-H-M-I1}  yields
\begin{equation*}
- \frac{T}{4\eps} \le - \frac{\kappa_\alpha}{a^{\frac{1}{\alpha}} \eps} + \ln \| u \|_{L^1(0, T)}  + \frac{T}{4 \eps} -  C \ln \eps +C.  
\end{equation*}
This implies
\begin{equation*}
\ln  \| u \|_{L^1(0, T)}   \ge \frac{\kappa_\alpha}{a^{\frac{1}{\alpha}} \eps} - \frac{T}{2 \eps}  + C \ln \eps - C. 
\end{equation*}
The conclusion follows.  
\end{proof}

Here is the second lemma, which is a consequence of \Cref{lem-H-M}.

\begin{lemma} \label{lem-HH} Let $\alpha > 1$, $a > 0$,  and let $(\lambda_k)_{k \ge 1}$ be a positive sequence. 
Assume that $\gamma > 0$ and $\Gamma_1 < + \infty$,  where 
\be
\gamma = \inf_{k \neq n} |\lambda_k - \lambda_n| \quad \mbox{ and } \quad \Gamma_1 = \sup_{k} \frac{|\lambda_k - a k^\alpha|}{k^{\alpha - 1}}.  
\ee
Let $0< T < m$ and  $u \in L^1(0, T)$. Define 
\be
\cG (z) = \int_0^T u(t) e^{- i z t} \, dt. 
\ee
Assume that $\cG(\lambda_1) = 1$ and $\cG (\lambda_k) = 0$ for $k \ge 2$. Then 
\begin{equation}
\| u \|_{L^1(0, T)}   \ge  c T^{C} e^{\frac{\rho}{T^{\tau}}} \mbox{ with } \tau = \frac{1}{\alpha - 1},  
\end{equation}
for some  positive constants $c$, $C$ depending only on $a, \, \alpha, \,  \gamma, \, \Gamma_1$, and $m$, and 
$$
\rho = 2^{\frac{1}{\beta}} \beta \left(  \frac{\kappa_\alpha}{a^{1/\alpha} (1 + \beta)}\right)^{1 + \beta},  
$$
with $\kappa_\alpha$ being defined by \eqref{def-ds}.  

\end{lemma}

\begin{proof} The proof is similar to the one of \Cref{lem-S} and the details are omitted. 
\end{proof}

\subsection{Lower bound on the cost of fast controls for a singular heat control system}

In this section, we deal with \eqref{sys-H-singular}. 

\begin{proposition}\label{pro-S-M} Let $0< \eps < m$ and $1/2 < s \le 1$. Denote  $\beta = 2s - 1$ and consider the control system, for $T>0$, 
\begin{equation}\label{pro-S-M-sys} \left\{
\begin{array}{c}
i y_t -  \eps^\beta (-\Delta)^s y + \frac{i}{2 \eps } y = 0   \mbox{ in } (0, T) \times (0, L), \\[6pt]
y(t, 0) = 0 \mbox{ in } (0, T), \\[6pt]
y (t, L) = u(t) \mbox{ in } (0, T), \\[6pt]
y_0 = \varphi_{1} \mbox{ in } (0, L). 
\end{array} \right. 
\end{equation} 
Assume that 
\be
y(T, \cdot) = 0 \mbox{ in } (0, L).   
\ee
Then 
\begin{equation}
\| u \|_{L^1(0, T)}   \ge c \eps^C e^{\frac{\theta_{2s}/\pi - T/2}{\eps}} 
\end{equation}
for some positive constants $c$, $C$ depending only on $s$ and $m$. Here $\theta_{2s}$ is defined by \eqref{def-cs}.  
\end{proposition}

The control system \eqref{pro-S-M-sys} is understood under the sense that $y \in C([0, T]; L^2(0, 1))$, and $y$ satisfies 
\be \label{pro-S-M-meaning}
i \frac{d}{dt} \langle y, \varphi \rangle - \eps^{\beta} \langle y, (-\Delta)^s \varphi \rangle +  \frac{i}{2 \eps } \langle y, \varphi \rangle = u(t) \varphi_x(t, 1) \mbox{ in } (0, T) \mbox{ for $\varphi = \varphi_n$ with $n \ge 1$}, 
\ee
in the distributional sense.

\begin{proposition}\label{pro-H-M} Let $0< \eps < m$ and $1/2 < s  \le 1$. Denote  $\beta = 2s  - 1$ and consider the control cost of the system, for $T>0$,  
\begin{equation} \label{pro-H-M-sys} \left\{
\begin{array}{c}
y_t +  \eps^\beta (-\Delta)^s y + \frac{1}{2 \eps } y = 0   \mbox{ in } (0, T) \times (0, L), \\[6pt]
y(t, 0) = 0 \mbox{ in } (0, T), \\[6pt]
y (t, L) = u(t) \mbox{ in } (0, T), \\[6pt]
y_0 = \varphi_{1} \mbox{ in } (0, L). 
\end{array} \right. 
\end{equation} 
Assume that 
\be
y(T, \cdot) = 0 \mbox{ in } (0, L).   
\ee
Then 
\begin{equation*}
\ln  \| u \|_{L^1(0, T)}   \ge c \eps^C e^{\frac{\kappa_{2s}/\pi - T/2}{\eps}} 
\end{equation*}
for some positive constants $c$, $C$ depending only on $s$ and $m$, where  $\kappa_{2s}$ being defined by \eqref{def-ds}.   
\end{proposition}

The control system \eqref{pro-H-M-sys} is understood under the sense that $y \in C([0, T]; L^2(0, 1))$, and $y$ satisfies  
\be \label{pro-H-M-meaning}
\frac{d}{dt} \langle y, \varphi \rangle + \eps^\beta \langle y, (-\Delta)^s \varphi \rangle + \frac{1}{2 \eps} \langle y, \varphi \rangle  = -  u(t)  \varphi_x(t, 1) \mbox{ in } (0, T) \mbox{ for $\varphi = \varphi_n$ with $n \ge 1$}, 
\ee
in the distributional sense.

\begin{proof} Define, for $k \ge 1$, 
\begin{equation*}
\Gamma_k :=   \eps^\beta \lambda_k +  \frac{1}{2 \eps}, 
\end{equation*}
where 
$$
\lambda_k = \gamma_k^{\alpha/2} \mbox{ with } \alpha = 2 s.  
$$

Consider \eqref{pro-H-M-meaning} with $\varphi = \varphi_k$. Multiplying this system with $e^{ \Gamma_k t}$ and integrating over $(0, T)$, we obtain
\be
 \langle y_0, \varphi_k \rangle = (-1)^k k \pi \int_0^T u(t) e^{\Gamma_k t} \, dt. 
\ee
Set 
$$
\cG(z) = \int_0^T u(t) e^{- i z t} \, dt
$$
Then 
$$
\cG(i \Gamma_1) = - \pi  \quad \mbox{ and } \cG(i \Gamma_k ) = 0 \mbox{ for } k \ge 2.  
$$
The conclusion now follows from \Cref{lem-H-M}. 
\end{proof}

\subsection{Proof of part $a)$ of \Cref{thm-H}}
Denote 
$$
\beta = 2s  - 1, \quad T_0 = \frac{2 c_s}{\pi (\beta + 1)}, \quad \eps^\beta T_0 = T, 
$$
and set 
\be
\ty(t, x) = y(\eps^\beta t, x) e^{-\frac{t}{2 \eps}},  \mbox{ and } \quad \tu(t) = u(\eps^\beta t) e^{-\frac{t}{2 \eps}}. 
\ee
We then have 
\begin{equation} \left\{
\begin{array}{c}
\ty_t +  \eps^\beta (-\Delta)^\alpha \ty + \frac{1}{2 \eps } \ty = 0   \mbox{ in } (0, T_0) \times (0, L), \\[6pt]
\ty(t, 0) = 0 \mbox{ in } (0, T_0), \\[6pt]
\ty (t, L) = \tu(t) \mbox{ in } (0, T_0), \\[6pt]
\ty_0 = \varphi_{1} \mbox{ in } (0, 1). 
\end{array} \right. 
\end{equation} 
Moreover, 
$$
\ty(T_0, \cdot) = y (T, \cdot) e^{-\frac{T_0}{2 \eps}} = 0. 
$$
By \Cref{pro-H-M}, we have 
$$
\|\tu \|_{L^1(0, T_0)} \ge \eps^C e^{\frac{\kappa_{2s}}{\pi \eps} - \frac{T_0}{2 \eps}}. 
$$
Since 
\be
\frac{\kappa_{2s}}{\pi \eps} - \frac{T_0}{2 \eps} = \frac{\kappa_{2s} T_0^{1/\beta}}{\pi T^{1/ \beta}} -  \frac{T_0^{1 + 1/ \beta}}{2 T^{1/\beta}} = \frac{1}{T^{\tau}} \left( \frac{\kappa_{2s}}{\pi} \right)^{1 + 1/ \beta} \left(\frac{2}{\beta + 1}\right)^{1/\beta}  \frac{\beta}{1 + \beta} = \frac{\mu_s}{T^\tau}
\ee
and  
\be
\|u \|_{L^1(0, T)} \ge \eps^{-\beta} \|\tu \|_{L^1(0, T_0)}, 
\ee
the conclusion follows. \qed

\begin{remark} \rm
One can obtain a proof of part $a)$ of \Cref{thm-H} using \Cref{lem-HH}. 
\end{remark}

\begin{remark} \rm
Similar ideas as given in the proof of part $a)$ of \Cref{thm-S,thm-H} were implemented by the author in \cite{Ng-S-Schrodinger} to deal with the bilinear control of the Schr\"odinger equation.
\end{remark}

\section{Technical lemmas} \label{sec-lemmas}

In this section, we present several results which are the main ingredient in the proof of the upper bound of the cost in \Cref{thm-S,thm-H}. 
The first one is the following.

\begin{lemma}  \label{lem-Phi} Let $\alpha > 1$, $a > 0$,  and let $(\lambda_k)_{k \ge 1}$ be a positive sequence. 
Assume that $\gamma > 0$ and $\Gamma_1, \Gamma_2 < + \infty$ where 
\be \label{lem-Phi-assumption}
\gamma = \inf_{k \neq n} |\lambda_k - \lambda_n|, \quad  \Gamma_1 = \sup_{k} \frac{|\lambda_k - a k^\alpha|}{k^{\alpha - 1}}, \quad \mbox{ and } \quad \Gamma_2 = \sup \frac{k^\alpha}{\lambda_k}.  
\ee
Set, for $n \in \N$,  
$$
\Phi_n (z) =\prod_{k \neq n} \left( 1 - \frac{z}{\lambda_k - \lambda_n} \right) \mbox{ for } z \in \mC. 
$$
We have 
\be \label{lem-Phi-cl2}
|\Phi_n (z)| \le C e^{c  |z|^{\frac{1}{\alpha}}}  \mbox{ for large $|z|$ with $z \in \mC$,}   
\ee
and
\be \label{lem-Phi-cl1}
|\Phi_n (- i x - \lambda_n) | \le C e^{ c |x|^{1/\alpha} + C \lambda_n^{1/\alpha}} \mbox{ for } x \in \mR, 
\ee
for some positive constants $c$ and  $C$ depending only on $a$, $\gamma$, $\Gamma_1$, and $\Gamma_2$. \end{lemma}

\begin{remark} \rm Thus the conclusion of \Cref{lem-Phi} holds if $(\lambda_k)_{k \ge 1}$ be a positive sequence such that 
\be
\inf_{k \neq n} |\lambda_k - \lambda_n| > 0   \quad \mbox{ and } \quad \lambda_k =  a k^\alpha + O(k^{\alpha -1}). 
\ee
since $\gamma>0$ and  $\Gamma_1, \Gamma_2 < + \infty$ in this case. 
\end{remark}

\begin{proof} We only prove \eqref{lem-Phi-cl2} for $|x|> 1$ and \eqref{lem-Phi-cl1} for $|z|> 1$. The other cases are simpler and omitted. In what follows in this proof, $C$ denotes a positive constant depending only on $a$, $\gamma$, $\Gamma_1$, and  $\Gamma_2$, and can change from one place to another.

\medskip 
We begin with \eqref{lem-Phi-cl2}.  Given $z \in \mC$ with $|z| > 1$, let $m \in \N$ be such that $m^\alpha \le |z| < (m+1)^\alpha$.

We first consider the case $\lambda_n <  |z|$. We have
\begin{multline}\label{lem-Phi-p2-1}
\ln |\Phi_n (z)| \le C \sum_{|\lambda_k-\lambda_n|  \ge C |z|} \frac{|z|}{|\lambda_k-\lambda_n|} + C \sum_{|\lambda_k-\lambda_n|  < C |z|} \ln \left(1 +  \frac{|z|}{|\lambda_k - \lambda_n|} \right) \\[6pt] 
\le \frac{C |z|}{|z|^{\frac{\alpha - 1}{\alpha}}} + \sum_{1 \le n^{\alpha-1} m \le C |z|}  \ln \left( 1+  \frac{C |z|}{n^{\alpha-1}m} \right) + C.   
\end{multline}

We have, with $a = C |z|/ n^{\alpha-1} \ge 1$,  
\begin{multline}\label{lem-Phi-p2-2}
\int_1^{a} \ln (1 + a/ t) \, dt = \int_1^{a} \Big( \ln (a + t) - \ln t \Big) \, dt = \int_{1 + a}^{2a} \ln t \, dt - \int_1^{a} \ln t \, dt \\[6pt]
=  \int_{a+1}^{2a} \ln t \, dt -\int_1^{a} \ln t \, dt  = 2 a \ln (2a) - (a+1) \ln (1 + a) - a \ln a - \int_{a+1}^{2a} 1 \, dt + \int_1^{a} 1 \, dt 
\le C a. 
\end{multline}
Combining \eqref{lem-Phi-p2-1} and \eqref{lem-Phi-p2-2} yields 
$$
\ln |\Phi_n (z)| \le C  |z|^{1/ \alpha} + C,  
$$
which implies \eqref{lem-Phi-cl2}.

We next deal with the case $|z| \le |\lambda_n|$. We have
\begin{multline}\label{lem-Phi-p2-3}
\ln |\Phi_n (z)| \le C \sum_{|\lambda_k|  \ge 2 |\lambda_n| } \frac{|z|}{|\lambda_k|} + C \sum_{|\lambda_k|  \le 2 |\lambda_n|} \ln \left(1 +  \frac{|z|}{|\lambda_k - \lambda_n|} \right) + C \\[6pt] 
\le \frac{C |z|}{n^{\alpha-1}} + \sum_{1 \le m \le C n}  \ln \left( 1+  \frac{C |z|}{n^{\alpha-1} m} \right) + C.   
\end{multline}
We have, with $a = C |z|/ n^{\alpha-1} $ and $\hat n = C n$,  
\begin{multline}\label{lem-Phi-p2-4}
\int_1^{\hat n} \ln (1 + a/ t) \, dt = \int_1^{\hat n} \Big( \ln (a + t) - \ln t \Big) \, dt = \int_{1 + a}^{\hat n+ a} \ln t \, dt - \int_1^{\hat n} \ln t \, dt \\[6pt]
=  \int_{\hat n}^{\hat n+ a} \ln t \, dt -\int_1^{1+ a} \ln t \, dt  = (\hat n+ a) \ln (\hat n+a) - \hat n \ln \hat n  - (a+1) \ln (1 + a)  - \int_{\hat n}^{\hat n+a} 1 \, dt + \int_1^{1+ a} 1 \, dt 
\le C a. 
\end{multline}
Combining \eqref{lem-Phi-p2-3} and \eqref{lem-Phi-p2-4} yields, since $|z| \le |\lambda_n|$,  
$$
\ln |\Phi_n (z)| \le C  |z|^{1/ \alpha} + C,  
$$
which implies \eqref{lem-Phi-cl2}.

We next deal with \eqref{lem-Phi-cl1}. We have 
$$
|\Phi_n (- i x - \lambda_n) | = \prod_{k \neq n} \left| \frac{i x + \lambda_k}{\lambda_k - \lambda_n} \right| = \frac{\prod_{k \neq n} \left| 1+ i x / \lambda_k \right|}{\prod_{k \neq n}  \left| 1- \lambda_n / \lambda_k \right|}  \mbox{ for } x \in \mR. 
$$
Given $x \in \mR$ with $|x| > 1$, let $m \in \N$ be such that $m^\alpha \le |x| < (m+1)^\alpha$. We then have, by \eqref{lem-Phi-assumption}, 
\begin{multline} \label{lem-Phi-p1-1}
\ln \left( \prod_{k \neq n} \left| 1+ i x / \lambda_k \right|^2 \right)  \le \sum_{k \ge 1} \ln |1 + x^2/ \lambda_k^2| = \sum_{\lambda_k > m^\alpha} \ln |1 + x^2/ \lambda_k^2| +  \sum_{\lambda_k \le m^\alpha} \ln |1 + x^2/ \lambda_k^2| \\[6pt]
\le C \sum_{\lambda_k > m^\alpha } x^2/ \lambda_k^2+  \ln \left| \prod_{\lambda_k \le  m^\alpha} (2 x^2/ \lambda_k^2) \right| \\[6pt]
\le C x^{2} \sum_{k > m} \frac{1}{k^{2 \alpha}} + \ln \Big(C^m x^{2m}/ (m!)^{2 \alpha} \Big) + C \le 
C x^{2}/(m+1)^{2 \alpha - 1} + \ln \Big( C^{m} x^{2m}/m^{2\alpha m} \Big) + C
\\[6pt]
\le C |x|^{\frac{1}{\alpha}} + \ln \Big( C^{m} x^{2m}/m^{2\alpha m} \Big) + C \le C |x|^{\frac{1}{\alpha}} + C. 
\end{multline}
Since, by \eqref{lem-product-p2} of \Cref{lem-product}  below,  
\be \label{lem-product-remark}
\frac{\prod_{k \neq n} \left| 1+ i x / \lambda_k \right|}{\prod_{k \neq n}  \left| 1- \lambda_n / \lambda_k \right|} \le C e^{C \lambda_n^{1/\alpha}}, 
\ee
Assertion \eqref{lem-Phi-cl1} now follows from  \eqref{lem-Phi-p1-1}. 
\end{proof}

In the proof of \Cref{lem-Phi}, we used the following result. 

\begin{lemma} \label{lem-product}
Let $\alpha > 1$, $a > 0$,  and let $(\lambda_k)_{k \ge 1}$ be a positive (strictly) increasing sequence. 
Assume that 
\be \label{lem-product-assumption}
0< \gamma = \inf_{k \neq n} |\lambda_k - \lambda_n|, \quad  + \infty > \Gamma_1 = \sup_{k} \frac{|\lambda_k - a k^\alpha|}{k^{\alpha - 1}}.  
\ee
Then 
\begin{equation}\label{lem-product-p1}
\prod_{j=1}^{\infty}
\left(1+\frac{\lambda_n}{\lambda_j}\right)
=
\exp\!\left[
\big(a^{-1/\alpha}P_\alpha + o(1)\big)\lambda_n^{1/\alpha}
\right],
\qquad n\to\infty,
\end{equation}
\begin{equation}\label{lem-product-p2}
\prod_{j=1; j \neq n}^{\infty}
\left|1-\frac{\lambda_n}{\lambda_j}\right|
= \exp\!\left[
\big(a^{-1/\alpha}Q_\alpha + o(1)\big)\lambda_n^{1/\alpha}
\right],
\qquad n\to\infty, 
\end{equation}
for some positive constant $k$ depending only on $\Gamma_1$,  where 
\begin{equation} \label{lem-product-PQ}
P_\alpha
=
\int_0^{\infty}
\frac{dv}{v^{1-1/\alpha}(v+1)}  \quad \mbox{ and } \quad 
Q_\alpha
=
\text{p.v.}\int_0^\infty
\frac{dv}{v^{1-1/\alpha}(1-v)}.
\end{equation}
\end{lemma}

\begin{proof} The proof is divided into two steps.

\medskip 
\noindent{\it Step 1:} Proof of \eqref{lem-product-p1}. Let $N(u)$ be the counting function for the sequence
$\{\lambda_n\}$, i.e.,  
\begin{equation}
N(u)=
\begin{cases}
0, & 0 \le u < \lambda_1,\\
n, & \lambda_n \le u < \lambda_{n+1} \mbox{ for } n \ge 1, 
\end{cases}
\end{equation}
It follows from \eqref{lem-product-assumption} that, for some $M > 0$, 
\begin{equation}\label{lem-product-Nu}
|N(u) -  a^{-1/\alpha}u^{1/\alpha}| \le M. 
\qquad u \ge 0. 
\end{equation}

By the properties of the Stieltjes integral, we have
\begin{align}
\log\!\left(
\prod_{j=1}^{\infty}
\left(1+\frac{\lambda_n}{\lambda_j}\right)
\right)
&=
\sum_{j=1}^{\infty}
\log\!\left(1+\frac{\lambda_n}{\lambda_j}\right) \notag\\
&=
\int_{\lambda_1^-}^{\infty}
\log\!\left(1+\frac{\lambda_n}{u}\right)\, dN(u)
=
\lambda_n
\int_{\lambda_1}^{\infty}
\frac{N(u)\,du}{u(u+\lambda_n)} .
\tag{4.5}
\end{align}
From \eqref{lem-product-Nu}, we have
\begin{equation} \label{lem-product-1-p1}
\left|
\int_{\lambda_1}^{\infty}
\frac{N(u)\,du}{u(u+\lambda_n)}
-
a^{-1/\alpha}
\int_{\lambda_1}^{\infty}
\frac{u^{1/\alpha}\,du}{u(u+\lambda_n)}
\right|
\le
M \int_{\lambda_1}^{\infty}
\frac{du}{u(u+\lambda_n)} = M\ln\!\left(\frac{\lambda_n}{\lambda_1}+1\right)
= O(\ln n). 
\end{equation}
On the other hand
\begin{multline}  \label{lem-product-1-p2}
a^{-1/\alpha}\lambda_n
\int_{\lambda_1}^{\infty}
\frac{u^{1/\alpha}\,du}{u(u+\lambda_n)}
=
a^{-1/\alpha}\lambda_n^{1/\alpha}
\int_{\lambda_1/\lambda_n}^{\infty}
\frac{v^{1/\alpha}\,dv}{v(v+1)}\\[6pt]
=
a^{-1/\alpha}\lambda_n^{1/\alpha}
\left[
\int_0^{\infty}
\frac{dv}{v^{1-1/\alpha}(v+1)}
+o(1)
\right],
\qquad \lambda_n\to\infty .
\end{multline}
Assertion \eqref{lem-product-p1} now follows from \eqref{lem-product-1-p1} and \eqref{lem-product-1-p2}. 

\medskip 

\medskip 
\noindent{\it Step 2:} Proof of \eqref{lem-product-p2}. Let $N_n(u)$ be  the counting function of the ``deleted'' sequence
$\{\lambda_j; j\neq n\}$
It follows from \eqref{lem-product-assumption} that, for some $M > 0$, 
\begin{equation}\label{lem-product-Nnu}
|N_n(u) -  a^{-1/\alpha}u^{1/\alpha}| \le M. 
\qquad u \ge 0. 
\end{equation}

We have, for $k \ge 1$,  
\begin{multline} \label{lem-product-2-p1}
\log\!\left|
\prod_{|j - n| \ge k}
\left(1-\frac{\lambda_n}{\lambda_j}\right) \right|
\\[6pt]
=
\int_{\lambda_1^-}^{\lambda_{n-k}}
\log\!\left(\frac{\lambda_n}{u}-1\right)\, dN_n(u)
+
\int_{\lambda_{n+k}^-}^{\infty}
\log\!\left(1-\frac{\lambda_n}{u}\right)\, dN_n(u). 
\end{multline}
By integration by parts, 
\begin{multline} \label{lem-product-2-p2}
RHS \mbox{ of } \eqref{lem-product-2-p1} = 
\left[
N_n(u)\log\!\left(\frac{\lambda_n}{u}-1\right)
\right]_{\lambda_1^-}^{\lambda_{n-k}}
+
\left[
N_n(u)\log\!\left(1-\frac{\lambda_n}{u}\right)
\right]_{\lambda_{n+k}^-}^{\infty}
\\[6pt]
+
\lambda_n
\int_{\lambda_1}^{\lambda_{n-k}}
\frac{N_n(u)\,du}{u(\lambda_n-u)}
+
\lambda_n
\int_{\lambda_{n+k}}^{\infty}
\frac{N_n(u)\,du}{u(u-\lambda_n)} .
\end{multline}

We treat the boundary terms first. Since 
\[
\lim_{u\to\infty}
N_n(u)\log\!\left(1-\frac{\lambda_n}{u}\right)=0,
\]
the contribution of the boundary terms in \eqref{lem-product-2-p2} is
\begin{multline}\label{lem-product-2-p3}
(n-k) \ln \left(\frac{\lambda_n}{\lambda_{n-k}}-1\right) - (n+k)\ln \left(1 - \frac{\lambda_n}{\lambda_{n+k}}\right) \\[6pt]
= n \ln \left( \frac{(\lambda_{n} - \lambda_{n-k}) \lambda_{n+k}}{(\lambda_{n+k} - \lambda_n) \lambda_{n-k}} \right) + k \ln \left( \frac{(\lambda_{n} - \lambda_{n-k})(\lambda_{n+k} - \lambda_{n})}{\lambda_{n+k} \lambda_{n-k}}\right)  \mathop{=}^{\eqref{lem-product-assumption}} O(n). 
\end{multline}
for a fixed $k$ sufficiently large by \eqref{lem-product-assumption} ($k$ is chosen independently of $n$ for large $n$). This will be assumed from later on. 

As in \eqref{lem-product-1-p1},  we have
\begin{multline} \label{lem-product-2-p4}
\lambda_n \int_{\lambda_1}^{\lambda_{n-k}}
\frac{N_n(u)\,du}{u(\lambda_n-u)}
+ \lambda_n \int_{\lambda_{n+k}}^{\infty}
\frac{N_n(u)\,du}{u(\lambda_n-u)} \\[6pt]
=  \lambda_n a^{-1/\alpha}
\int_{\lambda_1}^{\lambda_{n-k}}
\frac{u^{1/\alpha}\,du}{u(\lambda_n-u)}
+ \lambda_n a^{-1/\alpha}
\int_{\lambda_{n+k}}^{\infty}
\frac{u^{1/\alpha}\,du}{u(\lambda_n-u)} + O (\ln n )
\end{multline}
and, as in \eqref{lem-product-1-p2}, we obtain 
\begin{equation}\label{lem-product-2-p5}
\lambda_n a^{-1/\alpha}
\int_{\lambda_1}^{\lambda_{n-k}}
\frac{u^{1/\alpha}\,du}{u(\lambda_n-u)}
+
\lambda_n a^{-1/\alpha}
\int_{\lambda_{n+k}}^{\infty}
\frac{u^{1/\alpha}\,du}{u(\lambda_n-u)}
=
\lambda_n^{1/\alpha} a^{-1/\alpha} \Big( Q_\alpha + o(1) \Big). 
\end{equation}
Assertion \eqref{lem-product-p2} now follows from \eqref{lem-product-2-p1}-\eqref{lem-product-2-p5}, and \eqref{lem-product-assumption}. 
\end{proof}

\begin{remark} \rm In the proof of \Cref{lem-Phi}, we used the inequality \eqref{lem-product-remark}. This inequality can be proved directly with a more elementary proof. Nevertheless, we present here its proof via \Cref{lem-product}, which are interesting in themselves and gives the optimal behavior of LHSs of \eqref{lem-product-p1} and \eqref{lem-product-p2}. \Cref{lem-product} is an improvement of \cite[Lemma 3.1]{FR71} where the same results are proved under the assumption that, for some $c_0 > 0$,  
$$
\lambda_k = a (k + c_0)^\alpha  + o(k^{\alpha - 1}) 
$$
instead of \eqref{lem-product-assumption}. The proof presented here improves some technique points in \cite{FR71} and has its root from there. 
\end{remark}

\begin{remark} \rm It is noted in  \cite{FR71} that 
$$
Q_2 = 0.  
$$
This is consistent with the results given in \cite{TT07} where a better estimate is given in this case. The proof in \cite{TT07} used the Euler formula of the function $\sin (\pi x)$. 
\end{remark}

Here is the last result  of this section. 

\begin{lemma}\label{lem-H} Let $\mu > 0$, $\nu \ge 1$, and $0 < \beta \le 1 $. Set  
$$
H(z) = \int_{-1}^1 \sigma(t) e^{- i \beta z t } \, dt \mbox{ for } z \in \mC, 
$$
where 
$$
\sigma (t) = e^{- \frac{\nu^\mu}{(1 - t)^\mu} -  \frac{\nu^\mu}{(1 + t)^\mu}} \mbox{ for } t \in (-1, 1).  
$$
We have 
\be \label{lem-H-cl1}
|H(x)| \le C e^{- c |\beta \nu x|^{\frac{\mu}{\mu + 1}}} \mbox{ for } x \in \mR, 
\ee
and 
\be \label{lem-H-cl2}
H(ix) \ge \frac{1}{4} e^{-2^{\mu+1} \nu^\mu}e^{\beta |x|/4} \mbox{ for } x \in \mR, 
\ee
for some positive constants $C, c$ depending only on $\mu$. 
\end{lemma}

\begin{proof} We follow the idea of Tenenbaum and Tucsnak in \cite{TT07}. We first prove \eqref{lem-H-cl1}.  By integration by parts, we have 
\be \label{lem-H-p1}
|H(x)| \le \frac{1}{|\beta x|^{j}} \| \sigma^{(j)} \|_{L^\infty(-1, 1)} \mbox{ for } j \ge 1. 
\ee
For $t \in (0, 1)$, set $\rho = (1 - t)$. By Cauchy's integral formula, we have, for $k > 1$ sufficiently large (the largeness of $k$ depending only on $\mu$),  
\be \label{lem-H-p2}
\sigma^{(j)}(t) = \frac{j!}{2 \pi i} \int_{|z-t| = \rho/k} \frac{\sigma (z)}{(z- t)^{j+1}} \, dz. 
\ee

We have 
\be \label{lem-H-p3}
|\sigma (z)| \ge e^{- \frac{C_\mu \nu^\mu}{\rho^\mu}} \mbox{for $z \in \mC$ such that $|z-t| = \rho/k$}. 
\ee
Using the Stirling formula and the fact $e^{x} \ge x^n/n!$ for $x \ge 0$ for $n \ge 1$,  we derive from \eqref{lem-H-p2} and \eqref{lem-H-p3} that 
\be \label{lem-H-p4}
|\sigma^{(j)}(t)| \le C \frac{j! e^{- \frac{C \nu^\mu}{\rho^\mu}}}{\rho^j} \le C^j \frac{j^j}{\rho^j}  \frac{(\rho^\mu)^{[j/\mu]+1} ([j/\mu]+1)!}{(\nu^\mu)^{[j/\mu]+1}}\\[6pt]
\le  \frac{(C j)^{j (1 + 1/\mu)}}{\nu^{j}} \mbox{ for } t \ge 0. 
\ee
Since $\sigma$ is en even function on $(-1, 1)$, \eqref{lem-H-p4} holds for $t \in (-1, 1)$.  From \eqref{lem-H-p1}, we obtain  
$$
|H(x)| \le \left( \frac{C j^{1 + 1/\mu}}{\beta x \nu} \right)^{j}. 
$$
This implies, by choosing $j = [\beta \nu x]+1$,  
$$
|H(x)| \le C e^{- C (\beta \nu x)^\frac{\mu}{\mu+1}}.
$$
Assertion \eqref{lem-H-cl1} is established. 

We next deal with \eqref{lem-H-cl2}. Since $\sigma$ is an even function, it suffices to prove \eqref{lem-H-cl2} for $x \ge 0$. We have, for $x \ge 0$,  
$$
H(ix) = \int_{-1}^1 \sigma(t) e^{\beta x t} \, dt \ge \int_{1/4}^{1/2} \sigma(t) e^{\beta x t} \, dt \ge \frac{1}{4} e^{-2^{\mu+1} \nu^\mu}e^{\beta x/4}, 
$$
which is \eqref{lem-H-cl2}. 

The proof is complete. 
\end{proof}

\section{Upper bound of the cost of the fast controls} \label{sec-upperbound}

This section consisting of two subsections is devoted to the proof of the upper bounds given in \Cref{thm-S,thm-H}. These are given in \Cref{sect-UB-S} and \Cref{sect-UB-H}, respectively. 

\subsection{Proof of part $b)$ of \Cref{thm-S}} \label{sect-UB-S}
Set, for $k \in \N$,  
$$
\varphi_k(x) = \sin (k \pi x) \mbox{ in } (0, 1) \quad \mbox{ and } \quad \Phi_k(t, x) = e^{ i \lambda_k t} \varphi_k (x) \mbox{ in } (0, T) \times (0, 1). 
$$
where $\lambda_k = \lambda_k^{\alpha/2}$ with 
\be
\alpha = 2s.
\ee
Consider \eqref{thm-S-meaning} with  $\varphi = \varphi_k$. Multiplying this system with $e^{  i \lambda_k t}$ and integrating over $(0, T)$, we obtain
\be
- i \langle y_0, \varphi_k \rangle = (-1)^k k \pi \int_0^T u(t) e^{i \lambda_k t} \, dt = (-1)^k k \pi e^{i \lambda_k T/2} \int_{-T/2}^{T/2} v(t) e^{i \lambda_k t} \, dt 
\ee
where 
$$
v(t) = u(t + T/2) \mbox{ for } t \in (-T/2, T/2). 
$$
We then search a function $V$ such that 
$$
V(\lambda_k ) = c_k, 
$$
where 
\be
V(z) = \int_{-T/2}^{T/2} v(t) e^{- i z t} \, dt
\ee
and 
\be
c_k =  - \frac{i \langle y_0, \varphi_k \rangle }{(-1)^k k \pi e^{i \lambda_k T/2}}. 
\ee

Set, for $n \ge 1$, 
$$
\Psi_n (z) = \Phi_n (z - \lambda_n) = \prod_{k \neq n} \left(\frac{1 - z/\lambda_k}{1 - \lambda_n/ \lambda_k} \right) \mbox{ for } z \in \mC, 
$$
and 
$$
g_n (z) = \Psi_n(-z) H (z + \lambda_n) \mbox{ for } z \in \mC,  
$$
where, with $\beta = T/2$ and $\nu \beta = \gamma > 0$ sufficiently large but fixed, 
$$
H(z) = \alpha_0 \int_{-1}^1 \sigma(t) e^{- i \beta z t } \, dt \mbox{ with } 
\sigma (t) = e^{- \frac{\nu^\mu}{(1 - t)^\mu} -  \frac{\nu^\mu}{(1 + t)^\mu}} \mbox{ for } t \in (-1, 1),
$$
where $\mu$ is chosen such that 
\be
\frac{\mu}{\mu+1} = \frac{1}{\alpha} \quad \mbox{i.e., } \quad \mu = \frac{1}{\alpha - 1},  
\ee
and $\alpha_0$ is chosen such that $H(0) = 1$.  

Then 
$$
\alpha_0 \le C e^{2 \nu^\mu} = C e^{2 \nu^{\frac{1}{\alpha - 1}}}
$$
and 
$$
g_n(-\lambda_k)  = \delta_{n, k}. 
$$

We have, by \eqref{lem-Phi-cl1} of \Cref{lem-Phi} and \Cref{lem-H},  
$$
|g_n(x)| \le C e^{2 \nu^{\frac{1}{\alpha - 1}}} e^{ c_1 |x+ \lambda_n|^{\frac{1}{\alpha}}} e^{-c_2 ( \gamma |x + \lambda_n|) ^\frac{1}{\alpha}}  \mbox{ for } x \in \mR. 
$$
This implies, since $\gamma$ is sufficiently large, 
$$
|g_n(x)| \le C e^{2 \nu^{\frac{1}{\alpha - 1}}} e^{ - |x+ \lambda_n|^{\frac{1}{\alpha}}}  \mbox{ for } x \in \mR. 
$$

We now define 
$$
V(z) = \sum_{n \ge 1} c_n g_n (z) \mbox{ for } z \in \mC. 
$$
Then, with $a_n = \langle y_0, \varphi_n \rangle$,  
\be
\| V \|_{L^2(\mR)} \le \sum_{n \ge 1} |c_n| \|g_n \|_{L^2(\mR)} \le  C e^{2 \nu^{\frac{1}{\alpha - 1}}} \sum_{n \ge 1} \frac{|a_n|}{|n|} \le C e^{2 \nu^{\frac{1}{\alpha - 1}}}  \| y_0\|_{L^2(0, 1)}. 
\ee
We also have, by \eqref{lem-Phi-cl2} of \Cref{lem-Phi}, 
$$
|V(z)| \le \sum_{n \ge 1} |c_n| |g_n(z)| \le C_{T, \alpha} \sum_{n \ge 1} |c_n| e^{|z| T/2} e^{c |z|^{\frac{1}{\alpha}}} \le C_{T, \alpha, \eps} e^{|z| (T/2 + \eps)}  \| y_0\|_{L^2(0, 1)} \mbox{ for } z \in \mC. 
$$
It follows from Paley-Wiener's theorem, see, e.g., \cite[19.3 Theorem]{Rudin-RC},  that there exists $v \in L^2(\mR)$ with $\supp v \subset [-T/2, T/2]$ such that 
\be
V(z) = \int_{-T/2}^{T/2} v(t) e^{- i z t} \, dt,
\ee
\be
\| v\|_{L^2(\mR)} \le C e^{ \frac{C}{T^{\frac{1}{\alpha - 1}}}} \| y_0\|_{L^2(0, 1)}, 
\ee
and 
\be
V(\lambda_k) = c_k. 
\ee
The conclusion follows by taking $u(t) = v(t - T/2)$ for $t \in (0, T)$. \qed


\subsection{Proof of part $b)$ of \Cref{thm-H}} \label{sect-UB-H}
 Set, for $k \in \N$ with $k \ge 1$,  
$$
\varphi_k(x) = \sin (k \pi x) \mbox{ in } (0, 1), 
$$
where $\lambda_k = \gamma_k^{\alpha/2}$ with 
\be
\alpha = 2s.
\ee
Consider \eqref{thm-S-meaning} with  $\varphi = \varphi_k$. Multiplying this system with $e^{ \lambda_k t}$ and integrating over $(0, T)$, we obtain
\be
 \langle y_0, \varphi_k \rangle = (-1)^k k \pi \int_0^T u(t) e^{\lambda_k t} \, dt = (-1)^k k \pi e^{\lambda_k T/2} \int_{-T/2}^{T/2} v(t) e^{\lambda_k t} \, dt 
\ee
where 
$$
v(t) = u(t + T/2) \mbox{ for } t \in (-T/2, T/2). 
$$
We then search a function $V$ such that 
$$
V(i \lambda_k ) = c_k, 
$$
where 
\be
V(z) = \int_{-T/2}^{T/2} v(t) e^{- i z t} \, dt
\ee
and 
\be
c_k =   \frac{\langle y_0, \varphi_k \rangle }{(-1)^k k \pi e^{- \lambda_k T/2}}. 
\ee

Set 
$$
\Psi_n (z) = \Phi_n (- iz - \lambda_n) = \prod_{k \neq n} \left(\frac{1 + i z/\lambda_k}{1 - \lambda_n/ \lambda_k} \right)
$$
and 
$$
g_n (z) = \Psi_n(z) \frac{H (z)}{H(i \lambda_n)},  
$$
where, with $\beta = T/2$ and $\nu \beta = \gamma > 0$ sufficiently large but fixed, 
$$
H(z) = \int_{-1}^1 \sigma(t) e^{- i \beta z t } \, dt  \mbox{ with } 
\sigma (t) = e^{- \frac{\nu^\mu}{(1 - t)^\mu} -  \frac{\nu^\mu}{(1 + t)^\mu}} \mbox{ for } t \in (-1, 1), 
$$
where $\mu$ is defined by  
\be
\frac{\mu}{\mu+1} = \frac{1}{\alpha} \quad \mbox{ i.e., } \quad \mu = \frac{1}{\alpha - 1}. 
\ee
Then 
$$
g_n(i \lambda_k)  = \delta_{n, k}. 
$$

We have, by by \eqref{lem-Phi-cl1} of \Cref{lem-Phi}  and \Cref{lem-H},  for $x \in \mR$, 
\be
|g_n(x)| \le C e^{2 \nu^{\frac{1}{\alpha - 1}}} e^{ c_1 |x|^{\frac{1}{\alpha}} + C n} e^{- c_2 ( \gamma |x|) ^\frac{1}{\alpha}}  e^{C \nu^\mu - \beta \lambda_n/ 4} \le C e^{C \nu^{\frac{1}{\alpha - 1}} + C n - c_2 ( \gamma |x|) ^\frac{1}{\alpha} - \frac{T}{8} \gamma_n}. 
\ee
Since, for $\eps > 0$, there exists $C_\eps > 0$ such that, by Young's inequality,    
\be
n \le \eps T^\alpha n^\alpha + \frac{C_\eps}{T^{\alpha - 1}}, 
\ee
it follows that, for $\gamma$ sufficiently large, 
\be
|g_n(x)| \le C e^{2 \nu^{\frac{1}{\alpha - 1}}} e^{ - |x|^{\frac{1}{\alpha}}}  \mbox{ for } x \in \mR. 
\ee

We now define 
$$
V(z) = \sum_{n \ge 1} c_n g_n (z) \mbox{ for } z \in \mC. 
$$
Then 
\be
\| V \|_{L^2(\mR)} \le \sum_{n \ge 1} |c_n| \|g_n \|_{L^2(\mR)} \le  C e^{2 \nu^{\frac{1}{\alpha - 1}}} \sum_{n \ge 1} \frac{|a_n|}{|n|} \le C e^{2 \nu^{\frac{1}{\alpha - 1}}}  \| y_0\|_{L^2(0, 1)}. 
\ee
We also have, by \eqref{lem-Phi-cl2} of \Cref{lem-Phi} and \Cref{lem-H}, 
$$
|V(z)| \le \sum_{n \ge 1} |c_n| |g_n(z)| \le C_{T, \alpha} \sum_{n \ge 1} |c_n| e^{|z| T/2} e^{c |z|^{\frac{1}{\alpha}}} e^{- \frac{T}{8} \lambda_n} \le C_{T, \alpha, \eps} e^{|z| (T/2 + \eps)}  \| y_0\|_{L^2(0, 1)} \mbox{ for } z \in \mC. 
$$
It follows from Paley-Wiener's theorem, see, e.g., \cite[19.3 Theorem]{Rudin-RC},  that there exists $v \in L^2(\mR)$ with $\supp v \subset [-T/2, T/2]$ such that 
\be
V(z) = \int_{-T/2}^{T/2} v(t) e^{- i z t} \, dt
\ee
\be
\| v\|_{L^2(\mR)} \le C e^{ \frac{C}{T^{\frac{1}{\alpha - 1}}}} \| y_0\|_{L^2(0, 1)}, 
\ee
and 
\be
V(\lambda_k) = c_k. 
\ee
The conclusion follows by taking $u(t) = v(t - T/2)$ for $t \in (0, T)$. \qed

\end{document}